\documentclass[12pt]{article}
\usepackage{amsfonts}
\usepackage{mathrsfs}
\usepackage{amsmath}
\usepackage{latexsym}
\usepackage{amssymb}
\usepackage{amsthm}
\usepackage{amscd}
\title{\bf Nichols algebras over classical Weyl groups
}
\author{  \small Zhengtang Tan $^{a}$, Weicai Wu $^{b}$,  Shouchuan Zhang $^{c}$ \\
\small $a$. College of Engineering and Design, Hunan Normal University, Changsha  410081,  P.R. China \\
\small $b$. College of Mathematics,  Hunan Institute of Science and Technology, Yueyang 414006,   P.R. China\\
\small $c$. Department  of Mathematics,   Hunan University,  Changsha  410082, P.R. China\\
\small Emails: weicaiwu@hnu.edu.cn (WCW) }

\date {}
\begin{document}
\newtheorem{Proposition}{Proposition}[section]
\newtheorem{Theorem}[Proposition]{Theorem}
\newtheorem{Definition}[Proposition]{Definition}
\newtheorem{Corollary}[Proposition]{Corollary}
\newtheorem{Lemma}[Proposition]{Lemma}
\newtheorem{Example}[Proposition]{Example}
\newtheorem{Remark}[Proposition]{Remark}
\newtheorem{Conjecture}[Proposition]{Conjecture}
\newtheorem{Problem}[Proposition]{\quad Problem}

\maketitle

\begin{abstract} We show that  except in several  cases conjugacy classes of classical Weyl groups $W(B_n)$ and $W(D_n)$ are of type {\rm D}. We prove that except in three  cases Nichols algebras of irreducible Yetter-Drinfeld ({\rm YD} in short )modules over the classical Weyl groups are infinite dimensional.



\end {abstract}


\section{Introduction}\label{introduction}

Nichols algebras play a fundamental role in the classification of finite-dimensional complex pointed Hopf algebras by means of the lifting method developed by Andruskiewitsch and Schneider \cite {AS02}. In this context, given a group $ G$, an important step to classify all finite-dimensional complex pointed
Hopf algebras $H$ with group-like $G(H) = G$ is to determine all the pairs $(\mathcal O, \rho )$ such that the associated Nichols algebra $  \mathfrak B(\mathcal O, \rho )$ is finite-dimensional. Here the pairs $(\mathcal O, \rho )$ are such that $\mathcal O$ runs over all conjugacy classes of $G$ and $\rho $ runs over all irreducible representations of the centralizer of $g$ in $ G$, with $g \in \mathcal O$ fixed. In general, this is a difficult task since Nichols algebra is defined by generators and relations.
In practice, it is often useful to discard those pairs such that $\dim  \mathfrak B (\mathcal O, \rho ) = \infty$. There are properties of the conjugacy class $\mathcal O$ that imply that $\dim  \mathfrak B (\mathcal O, \rho )  = \infty$ for any $\rho $, one of which is the property of being of type {\rm D}. This is useful since it reduces the computations to operations inside the group and avoids hard calculations of generators and relations of the corresponding Nichols algebra.

This work contributes to the classification of finite-dimensional Hopf algebras over an algebraically closed field  of characteristic $0$. This problem was
posed by Kaplansky in 1975. The lifting method of Andruskiewitsch and Schneider  describes a way to classify finite-dimensional complex Hopf algebras. In \cite {AS10} they obtained the classification over the finite abelian groups whose order is relative prime with 210.
Nichols algebras of braided vector spaces  $(\mathbb CX,  c_q)$,  where $X$ is a
rack and $q$ is a 2-cocycle in $X$, were studied in \cite {AG03}.  It was shown
\cite {AFGV08,  AFZ09,  AZ07} that Nichols algebras  $\mathfrak B({\mathcal O}_
\sigma, \rho)$ over symmetry groups  are infinite dimensional,  except in a
number of remarkable cases corresponding to ${\mathcal O}_\sigma$. Two of the
present authors \cite{ZZ12} showed  that except in three cases Nichols algebras
of irreducible Yetter-Drinfeld (YD) modules over classical Weyl groups $A
\rtimes \mathbb S_n$ supported by $\mathbb S_n$ are infinite dimensional.
However, the classification has not been completed for Nichols algebras over
general classic Weyl groups $W(B_n)$ and $W(D_n)$.

Note that $\mathbb Z_2^n \rtimes \mathbb S_n$
is isomorphic to Weyl groups $W(B_n)$ and $W(C_n)$ of $B_n$ and
$C_n$ for $n>2$. If $A= \{a \in
     \mathbb Z_2^n \mid \ a = ({a_1}, {a_2}, \cdots, {a_n}) $  $\hbox {
with   } a_1 +a_2 + \cdots + a_n \equiv 0  \ ({\rm mod } \ 2) \}$,
then $A \rtimes \mathbb S_n$ is isomorphic to Weyl group
$W(D_n)$ of $D_n$ for $n>3$. Obviously, when $A= \{a \in
     {\mathbb Z}_2^n \mid \ a = ({a_1}, {a_2}, \cdots, {a_n}) $  $\hbox {
with  all } a_i \equiv 0 \ ({\rm mod } \ 2)\}$, $A \rtimes \mathbb S_n$ is isomorphic to
Weyl group of $A_{n-1}$ for $n>1$. Note that $\mathbb{S}_n$ acts on
$A$ as follows: for any $a\in A$ with $a = ({a_1}, {a_2},
\cdots, {a_n})$ and $\sigma, \tau \in \mathbb{S}_n$,
$\sigma \cdot a :=( {a_{\sigma^{-1}(1)}}, {a_{\sigma^{-1}(2)}}, \cdots,{a_{\sigma^{-1}(n)}} ).$

It is clear that
 \begin {eqnarray*}
&& (a, \sigma )^{- 1} = (- (a _{\sigma (1)}, a_{\sigma (2)}, \cdots  , a_{\sigma (n)}), \sigma^{ - 1}) = (- \sigma ^{ - 1}(a), \sigma ^{- 1}), \\
&&(b, \tau  )(a, \sigma )(b, \tau  )^{- 1} = (b + \tau  (a) -  \tau  \sigma \tau ^{ - 1}(b), \tau  \sigma \tau ^{ - 1}).\\
&& a \sigma
=(a, \sigma).
\end  {eqnarray*}
Without specification,  $K_n := \{ a\in \mathbb Z_n \mid a_1 +a_2 +\cdots + a_n =0\}$ and $ K_n \rtimes  \mathbb{S}_n$  is a subgroup of $\mathbb Z_2^n \rtimes \mathbb S_n$. Let $\mathbb W_n$ denote $  K_n \rtimes  \mathbb{S}_n$  or $\mathbb Z_2^n \rtimes \mathbb S_n$ throughout this paper. When we say  type
{\rm D} we mean a rack or conjugacy class of type
{\rm D} and not a Weyl group of type
${\rm D}_n$.

In this paper we prove  that except in several  cases conjugacy classes of classical Weyl groups $\mathbb W_n$ are of type
{\rm D}, and  except in three  cases Nichols algebras of irreducible
YD modules over the classical Weyl groups are infinite dimensional.

The work is organized as follows. In Section \ref{notation} we provide some preliminaries and set our notations.
In Section \ref{extension}  we determine when the conjugacy classes of  juxtapositions of two elements
are of type {\rm D}. In  Section \ref{rack}  we examine the rack and
conjugacy classes of $\mathbb W_n$.
In  Section \ref{nichols algebra} we classify Nichols algebras of irreducible {\rm YD} modules over the classical Weyl groups.

\section*{Preliminaries and conventions}\label{notation}

Let $\mathbf{ k}$ be an algebraically closed field of characteristic zero and $(V, C)$
 a braided vector space, i.e.   $V$ is
 a  vector space over $\mathbf{ k}$ and $C \in {\rm Aut} (V \otimes V)$ is a solution
 of the braid equation $   ({\rm id} \otimes C_{V,V}) (C_{V,V}\otimes {\rm id})({\rm id} \otimes C_{V,V}) = (C_{V,V}\otimes {\rm id})({\rm id} \otimes C_{V,V})( C_{V,V}\otimes{\rm id})$.
 $C_{i, i+1}:=  {\rm id} ^{\otimes i-1}     \otimes C_{V, V} \otimes {\rm id} ^{\otimes n-i-1}$.
Let $S_{m}\in {\rm End}_{\mathbf{ k} }(V^{\otimes m})$
and $S_{1, j}\in {\rm End}_{\mathbf{ k}}(V^{\otimes (j+1)})$ denote the maps $S_{m}=\prod \limits_{j=1}^{m-1}({\rm id}^{\bigotimes m-j-1}  \bigotimes S_{1, j})$ ,
 $S_{1, j}={\rm id}+C_{12}^{-1}+C_{12}^{-1}C_{23}^{-1}+\cdots+C_{12}^{-1}C_{23}^{-1}\cdots C_{j, j+1}^{-1}$
(in leg notation) for $m\geq 2$ and $j\in \mathbb N$. Then the subspace
$S=\bigoplus \limits _{m=2}^{\infty}  {\rm ker} S_{m}$ of the tensor $T(V)=\bigoplus \limits _{m=0}^{\infty}V^{\otimes m}$
is a two-sided ideal,  and algebra $\mathfrak B(V)=T(V)/S$ is termed the Nichols algebra associated to $(V, C)$.

For  $s\in G$ and  $(\rho,  V) \in  \widehat {G^s}$,  here is a
precise description of the {\rm YD} module $M({\mathcal O}_s,
\rho)$,  introduced in \cite {Gr00}. Let $t_1 = s, t_2,   \cdots,
t_{m}$ be a numeration of ${\mathcal O}_s$,  which is a conjugacy
class containing $s$,   and let $g_i\in G$ such that $g_i \rhd s :=
g_i s g_i^{-1} = t_i$ for all $1\le i \le m$. Then $M({\mathcal
O}_s,  \rho) = \oplus_{1\le i \le m}g_i\otimes V$. Let $g_iv :=
g_i\otimes v \in M({\mathcal O}_s, \rho)$,  $1\le i \le m$,  $v\in V$.
If $v\in V$ and $1\le i \le m $,  then the action of $h\in G$ and the
coaction are given by
\begin {eqnarray} \label {e0.11}
\delta(g_iv) = t_i\otimes g_iv,  \qquad h\cdot (g_iv) =
g_j(\nu _i(h)\cdot v),  \end {eqnarray}
 where $hg_i = g_j\nu _i(h)$,  for
unique  $1\le j \le m$ and $\nu _i(h)\in G^s$. The explicit formula for
the braiding is then given by
\begin{equation} \label{yd-braiding}
C(g_iv\otimes g_jw) = t_i\cdot(g_jw)\otimes g_iv =
g_{j'}(\nu _j (t_i)\cdot w)\otimes g_iv\end{equation} for any $1\le i, j\le
m$,  $v, w\in V$,  where $t_ig_j = g_{j'}\nu _j(t_i)$ for unique $j'$,  $1\le
j' \le m$ and $\nu _j(t_i) \in G^s$. Let $\mathfrak{B} ({\mathcal O}_s,
\rho )$ denote $\mathfrak{B} (M ({\mathcal O}_s,  \rho ))$.
$M({\mathcal O}_s,  \rho )$ is a simple {\rm YD} module (see \cite
 {DPR91, Ci97,    AZ07}).

We briefly recall the definition and main properties of racks;
see \cite {AG03} for details, more information and bibliographical references.
A rack is a pair $(X, \rhd )$,  where $ X$
is a non-empty set and  $\rhd :  X \times X \rightarrow  X$ is an operation such that
  $ x\rhd x = x$, $x\rhd (y \rhd z) = (x\rhd y)\rhd (x\rhd z)$  and $\phi _x$
is invertible for any $x, y, z
\in  X$, where $\phi _x$ is a map from $X$ to $X$ sending $y$ to $x\rhd y$  for any $x, y\in X.$
For example, $({\mathcal O}_s^G, \rhd )$ is a rack  with $x\rhd y:= x y x^{-1}.$

 If $R$ and $S$ are two subracks of $X$  with $R\cup S = X$, $R\cap S =\emptyset$,  $x\rhd y \in S$, $y\rhd x \in R$,  for any $x\in R, y\in S,$ then $R\cup S$ is called a decomposition of subracks of $X$. Furthermore, if there exist $a\in R$, $b\in S$ such that ${\rm sq} (a, b) := a\rhd (b \rhd (a \rhd b)) \not= b$, then $X$ is called to be of type {\rm D}. Notice that if a rack $Y$ contains a subrack  $X$ of type ${\rm D}$, then $Y$  is also of  type ${\rm D}$ ( see \cite {AFGV08}).

\section{Conjugacy classes of juxtapositions}\label {extension}
In this section we determine when the conjugacy classes of  juxtapositions of two elements are of type {\rm D}.

\begin {Lemma}\label {1.1} Let $G= \mathbb W_n$ and $(a, \tau), (b, \mu)\in G$. Let $(c, \lambda)$ denote $(a, \tau) \rhd ((b, \mu) \rhd((a, \tau)\rhd (b, \mu)))$.

{\rm (i) } Then  \begin {eqnarray}\label {e2.1.1}  c &=&(a+\tau \cdot [ b+\mu  \cdot(a+\tau \cdot b+(\tau \rhd \mu )\cdot a)
+ (\mu \rhd (\tau \rhd \mu)) \cdot b ] \nonumber \\
& +& (\tau \rhd (\mu \rhd (\tau \rhd \mu))) \cdot a.
\end {eqnarray}

 {\rm (ii)} If  $\tau $ and $\mu$ are commutative,  then
\begin {eqnarray}\label {e2.1.2}  c=a+\tau \mu \cdot a+\tau \mu^2 \cdot a+\mu \cdot a+\tau \cdot b+\tau ^2 \mu \cdot b+\tau \mu \cdot b.\end {eqnarray}

{\rm (iii)} If $\tau $ and $\mu$ are commutative,  then ${\rm sq }((a, \tau),  (b, \mu)) = (b, \mu) $ if and only if
\begin {eqnarray}\label {e2.1.3}  a+\tau \mu \cdot a+\tau \mu^2 \cdot a+\mu\cdot a=b+\tau \cdot b+\tau ^2 \mu \cdot b+\tau \mu \cdot b.\end {eqnarray}

{\rm (iv)} If $\tau $ and $\mu$ are commutative with  $(b, \mu)=\xi\rhd (a, \tau )$ and $\xi \cdot a=a$,  then $c=a+\tau \mu^2 \cdot a+\mu \cdot a+\tau \cdot a+\tau ^2 \mu \cdot a$.

{\rm (v)} If  $\tau $ and $\xi$ are commutative with $\tau ^2=1$ and $(b, \mu)=\xi \rhd (a, \tau )$ and $\xi\cdot a=a$,  then $c=a$.
\end {Lemma}
\noindent {\it Proof.}  It is clear. \hfill $\Box$

\vskip.1in

If the lengths of independent  cycles of $ \pi$ and $\tau$ are different,  then they are called mutually orthogonal,  written as $\pi \bot \tau. $
If the lengths of independent sign cycles of $(a, \pi)$ and $(b, \tau)$ are different,  then they are called mutually orthogonal,  written as $a\pi \bot b\tau. $ Obviously, $a\pi \bot b\tau $ if and only if $\pi \bot \tau. $ If ${\rm ord (\pi)}$ and ${\rm ord (\tau)}$ are coprime  and one of the two elements does not have any fixed point, then  $\pi \bot \tau. $ Notice
that if $\pi \bot \tau $, then  $\mathbb S_{m+n}^{\pi\#\tau}=\mathbb S_n^{\pi}\#\mathbb S_m^{\tau}$.

For any $a\pi\in \mathbb W_n $ and $b\tau\in \mathbb W_m$,  define $a\pi \# b\tau \in \mathbb W_{m+n}$ as follows:

$(a\# b)_{i} :=\left
\{\begin
{array} {ll}  a_i &\hbox {when } i\le n\\
b_{i-n} &\hbox {when } i>n
\end {array} \right.,  $  $  (\pi \# \tau ) (i):=\left
\{\begin
{array} {ll} \pi (i) &\hbox {when } i\le n\\
\tau (i-n)+n &\hbox {when } i>n
\end {array} \right. $
\hbox { and } $ a\pi\#b \tau : = (a\# b,   \pi \# \tau ) $.
Obviously $a\pi\#b \tau \in \mathbb W_{m+n}$ and it is called a juxtaposition of $a\pi$ and $b \tau$  (see \cite {AZ07}). Let $\overrightarrow{\nu_{n,  m}}$ be a map from $\mathbb W_n$ to $\mathbb W_{m+n}$ by sending $a\pi$ to $\overrightarrow{\nu_{n,  m}}(a\pi) :=a\pi \# 1_{\mathbb W_m}$; let $\overleftarrow {\nu_{n,  m}}$ be a map from $\mathbb W_m$ to $\mathbb W_{m+n}$ by sending $b\tau$ to $\overleftarrow{\nu_{n,  m}}(b\tau) :=1_{\mathbb W_n} \# b \tau$.

\begin {Lemma}\label {2.4} Assume  $a\pi \bot b\tau$ with $a\pi,  a'\pi'\in \mathbb W_n $ and $b\tau,  b'\tau'\in \mathbb W_m$. Then

 {\rm (i)} $(a\pi\#b\tau)(a'\pi'\#b'\tau')=(a\pi a'\pi'\# b\tau b'\tau')$.

 {\rm (ii) } $a\pi\#b\tau=\overrightarrow{\nu_{n, m}}(a\pi)\overleftarrow{\nu_{n, m}}(b\tau)=\overleftarrow{\nu_{n,  m}}(b\tau)
\overrightarrow{\nu_{n, m}}(a\pi)$.

{\rm (iii)} $\mathbb W_{m+n}^{a\pi\#b\tau}=\mathbb W_n^{a\pi}\#\mathbb W_m^{b\tau}=\overrightarrow{\nu_{n, m}}(\mathbb W_{n} ^{a\pi})\overleftarrow{\nu_{n, m}}(\mathbb W_m^{b\tau})$ as directed products.

{\rm (iv)} For any $\rho\in\widehat{\mathbb W_{m+n}^{a\pi\# b\tau}}$,  there exist $\rho_1 \in\widehat{\mathbb W_{n}^{a\pi}}$,  $\rho _2   \in\widehat{\mathbb W_{m}^{b\tau}}$  such that $\rho=\rho_1\otimes\rho_2 $.

{\rm (v)} $(a\pi\# b\tau)\rhd(a'\pi'\# b'\tau')=(a\pi \rhd a'\pi')\#(b\tau\rhd b'\tau')$.

{\rm (vi)} $\mathcal O_{a\pi\# b\tau}^{\mathbb W_{m+n}}=\mathcal O_{a\pi}^{\mathbb W_n}\#\mathcal O_{b\tau}^{\mathbb W_m}$.
\end {Lemma}
\noindent {\it Proof.} {\rm (i)},  {\rm (ii)} and {\rm (v)} are clear.

{\rm (iii)} By \cite [Section 2.2] {AZ07},  $\mathbb S_{m+n}^{\pi\#\tau}=\mathbb S_n^{\pi}\#\mathbb S_m^{\tau}.$
Obviously,  $\mathbb W_n^{a\pi}\#\mathbb W_m^{b\tau}\subseteq \mathbb W_{m+n} ^{a\pi \# b\tau} $. For any
$c\xi\in\mathbb W_{m+n}^{a\pi\# b\tau}, $ then $\xi\in\mathbb S_{m+n}^{\pi \#\tau}$ and there exist $\mu\in\mathbb S_n^\pi$ and $\lambda\in\mathbb S_m ^\tau $ such that $\xi=\mu \# \lambda$. Consequently,  $c\xi=d\mu\#f\lambda.$ Considering $c\xi(a\pi \# b \tau)=(a\pi \# b\tau)c\xi$,  we have $d\mu \in \mathbb W_n^{a\pi}$ and
 $f\lambda  \in \mathbb W_m^{b\tau}$. This completes the proof.

{\rm (iv)} It follows from {\rm (iii)}.

(vi) By {\rm (v)},  $  \mathcal O_{a\pi}^{\mathbb W_{n}}\#\mathcal O_{ b\tau}^{\mathbb W_{m}}\subseteq\mathcal O_{a\pi\# b
\tau}^{\mathbb W_{m+n}}. $  Consequently {\rm (vi)} follows from {\rm (iii)}. \hfill $\Box$

\begin{Remark} {\rm (i),  (ii)} and {\rm (v)} above still hold when $a\pi$ and $b\tau$ are not mutually orthogonal. \end{Remark}

\begin {Theorem}\label {2.5} If $\mathcal O_{a\tau}$ is of type {\rm D},  then $\mathcal O_{a\tau \# b\mu}$ is also of type {\rm D}.
\end {Theorem}
\noindent {\it Proof.} Let $X= R\cup S$ be a subrack decomposition of $\mathcal O_{a\tau}$ and of type {\rm D}. It is clear that
$X\# b\mu=R\# b\mu\cup S \# b\mu$ is a subrack decomposition of $ \mathcal O_{a\tau\# b\mu}$ and of type {\rm D}. \hfill $\Box$

\begin {Lemma}\label {2.7} Assume $a\pi \bot b\tau$. Let $\rho=\rho_1\otimes \rho_2 \in \widehat {\mathbb W_{m+n}^{a\pi \# b\tau}}$,  $ \rho_1\in \widehat{\mathbb W_{n}^{a\pi }} $,  $\rho_2 \in \widehat {\mathbb W_{m}^{b\tau}}$,  $a\pi\in \mathbb W_n$ and $b\tau\in \mathbb W_m$ with  $q_{a\pi, a\pi} {\rm id} = \rho_1 (a\pi)$ and $q_{b\tau, b\tau } {\rm id}=\rho_2 (b\tau)$. If $\dim \mathfrak{B} ({\mathcal O}_{a \pi\# b \tau }^ {\mathbb W_{n+m} }, \rho)<\infty $,  then

{\rm (i)} $M( {\mathcal O}_{a \pi }^{\mathbb W_n},  \rho_1  )$ is isomorphic to a {\rm YD} submodule of  $M ({\mathcal O}_{a\pi\# b\tau }^  {\mathbb W_{n+m} },  \rho)$ over $\mathbb W_n$ when $q_{b\tau, b\tau}=1$; hence $\dim \mathfrak{B} ({\mathcal O}_{a \pi}^  {\mathbb W_{n} }, \rho_1) <\infty$.

{\rm (ii)} $q_{a\pi, a \pi} q_{b\tau, b\tau}=-1$.

{\rm (iii)} $q_{b\tau, b\tau}=1$ and $q_{a\pi, a\pi}=-1$ when ${\rm ord}(b\tau)\le 2$ and ${\rm ord}(q_{a\pi, a\pi})\not=1$.

{\rm (iv)} $q_{b\tau, b\tau} =1$ and $q_{a\pi, a\pi}=-1$ when ${\rm ord}(a\pi)$ and ${\rm ord }(b\tau)$ are coprime and ${\rm ord}(b\tau)$ is odd.
\end {Lemma}
\noindent {\it Proof.} {\rm (i)} The proof is similar to the one in \cite [Section 2.2] {AZ07}. Indeed,  Let $t_1 = \pi, t_2,   \cdots,
t_{\mid {\mathcal O}_\pi\mid }$ be a numeration of ${\mathcal O}_\pi$,    and let $g_i\in G$ such that $g_i \rhd \pi :=
g_i \pi g_i^{-1} = t_i$ for all $1\le i \le \mid {\mathcal O}_\pi\mid $. Let $s_1 = \tau, s_2,   \cdots,
s_{ \mid {\mathcal O}_\tau \mid}$ be a numeration of ${\mathcal O}_\tau$,    and let $h_i\in G$ such that $h_i \rhd \tau :=
h_i \tau  h_i^{-1} = s_i$ for all $1\le i \le  \mid  {\mathcal O}_\tau \mid$. Thus  $\{ t_i \# s_j \mid 1 \le i \le \mid {\mathcal O}_\pi\mid, 1\le  j \ \le \   \mid {\mathcal O}_\tau \mid  \} = {\mathcal O}_{\pi \#\tau}$.  Let  $V$ and $W$ be  representation spaces of $\rho_1$ and $\rho_2$, respectively. Let  $0\not= w_0\in W$. Define a map $\psi $
from  $M({\mathcal O}_\pi^{{\mathbb W}_n },
\rho_1)$ to $M({\mathcal O}_ {\pi \# \tau}^{{\mathbb W}_m  \# {\mathbb W}_n},
\rho_1 \otimes \rho_2)$  by sending $ g_i v$ to $(g_i \# h_1) v\otimes w_0$ for any $v\in V, $ $1\le i \le\mid {\mathcal O}_\pi\mid.$  It is clear that $\psi$ is injective. Now we show that $\psi$ is a homomorphism of braided vector spaces. For any $v, v'\in V,$ we have that
\begin {eqnarray*}  (\psi \otimes \psi) ( C (g_iv  \otimes g_jv')) &=&
(\psi \otimes \psi )(  g_{i\rhd j} \rho_1 (\nu _j(t_i))(v') \otimes g_i v  ) \\
&=& (g_{i\rhd j} \# h_1 \rho_1 (\nu  _j(t_i)(v') \otimes w_0) \otimes (g_i \# h_1 v \otimes w_0) \ \ \ \hbox {and }
\end {eqnarray*}
\begin {eqnarray*} C (\psi (g_iv) \otimes \psi (g_jv') &=&
C( (g_i \# h_1 v\otimes w_0) \otimes  (g_j \# h_1 v'\otimes w_0)   )\\
&=&  ( g_{i\rhd j} \# h_1( \rho_1 \otimes \rho_2 ) (  \nu _{j, 1} (t_i\# s_1) ) (v'\otimes w_0) ) \otimes (g_i \#h_1 v \otimes w_0).
\end {eqnarray*}
It is clear that $\nu _  {j, 1} (t_i \# h_1) = \nu _{ j} (t_i) \# h_1$. Therefore $\psi$
is a homomorphism of braided vector spaces.

{\rm (ii)} and  {\rm (iii)} are clear.

{\rm (iv)} Obviously,  ${\rm ord}(q_{a\pi, a\pi}) \mid {\rm ord}(a\pi)$  and ${\rm ord}(q_{b\tau, b\tau }) \mid {\rm ord}(b\tau)$. Therefore,  ${\rm ord}(q_{b\tau, b\tau})$ is odd.  The least  common factor $({\rm ord} (q_{a\pi, a\pi}),  {\rm ord}(q_{b\tau, b\tau}))=1$ since $({\rm ord}{(a\pi}), {\rm ord}({b\tau}))=1$. By Part {\rm (ii)},   ${\rm ord}(q_{a\pi, a\pi}){\rm ord}(q_{b\tau, b\tau})=2$. Consequently,  ${\rm ord}(q_{a\pi, a\pi})=2$ and ${\rm ord}(q_{b\tau, b\tau})=1$.
\hfill $\Box$

\section{Conjugacy classes of $\mathbb W_n$}\label {rack}
In this section we  prove that except in several  cases conjugacy classes of classical Weyl groups $\mathbb W_n$ are of type {\rm D}.

\begin {Lemma}\label {1.3} Let $p$ be odd with $p\ge 5$ and $a\tau \in \mathbb W_n$ with $\tau = (1\ 2\cdots \ p)$. Then $\mathcal{O}_{a\tau}^{\mathbb W_n}$ is of type {\rm D}.
\end {Lemma}
\noindent {\it Proof.} {\rm (i)} Assume that $a\tau$ is a negative cycle (defined in \cite [Appendix ]{ZZ12}) and $b=(1, 0, \cdots, 0, a_{p+1}, a_{p+2}, \cdots, a_n)$. Thus   $a \tau $ and $b \tau$ are conjugate. We assume $a =  (1, 1, \cdots, 1, $ $ a_{p+1}, a_{p+2}, \cdots, a_n) $ without lost generality. Obviously,  the right hand side of (\ref {e2.1.3}) is non-vanishing for $\mu=\tau ^2$,  i.e. ${\rm sq}(a\tau, b\tau^2) \not= b\tau^2.$ Let $R:=\mathbb Z_2^n \rtimes\tau \cap  \mathcal O _{a \tau}
^{\mathbb W _n} $   and $S:=\mathbb Z_2^n \rtimes \tau^2 \cap\mathcal O _{a \tau} ^{\mathbb W _n}$. It is clear that $R \cup
S$ is a subrack decomposition of   $\mathcal O_{a\tau }^{\mathbb W_n}$, hence it is of type {\rm D}; notice that $\tau$ and $\tau ^2$ are conjugate in $\mathbb S_n.$

{\rm (ii)} Assume that $a\tau$ is a positive cycle. Let $b=(1, 0, 0, 1, 0, a_{p+1}, a_{p+2}, \cdots, a_n)$ when $p=5;$ $b=(1, 1, 0, \cdots, 0, a_{p+1}, a_{p+2}, \cdots, a_n)$ when $p>5.$ Thus $0\tau$,  $a \tau $ and $b \tau$ are conjugate. We assume $a=0$ without lost generality. It is clear that the right hand side of (\ref
{e2.1.3}) is not equal to $0$. Consequently,  $R \cup S$ is a subrack decomposition and is of type {\rm D} as Part {\rm (i)}.

\hfill $\Box$

\begin {Lemma}\label {1.3'} If $\sigma$  is of type ($3^2$),  then $\mathcal O_{a\sigma} ^{\mathbb W_6}$ is of type ${\rm D}$ for all $a \in \mathbb Z_2^6.$
\end {Lemma}
\noindent {\it Proof.} Let $\pi = (1 \ 2\ 3)$, $\xi = (4\ 5\ 6)$, $\tau =(1\ 2\ 3)(4\ 5\ 6)$ and  $\mu = (1\ 2\ 3)^2(4\ 5\ 6)$. We have that
the $4$-th,  $5$-th and $6$-th components of (\ref {e2.1.3}) are
\begin {eqnarray}\label {e1.3'}  (a_6+a_5,  a_4+a_6,  a_5+a_4) =  (b_6+b_5,  b_4+b_6,  b_5+b_4). \end {eqnarray}
By (\ref {e2.1.3}),
{\rm (i)} If $a = (1, 1, \cdots, 1)$ and $b=(1, 0, 0, 1, 0, 0)$,  then (\ref {e1.3'}) does not hold.

{\rm (ii)} If $a=0$ and $b=(0, 0, 0, 1, 1, 0)$,  then (\ref {e1.3'}) does not hold.

{\rm (iii)} If $a=(1, 0, 0, 0, 0, 0)$ and $b=(1, 0, 0, 1, 1, 0)$,  then (\ref {e1.3'}) does not hold.

{\rm (iv) } If $a=(0, 0, 0, 1, 0, 0)$ and $b=(0, 0, 0, 0, 1, 0)$,  then (\ref {e1.3'}) does not hold.
Let $R := \mathbb Z_2^6 \rtimes \tau \cap \mathcal O_{a\sigma}^{\mathbb W_6}$ and $S := \mathbb Z_2^6\rtimes \mu\cap \mathcal O_{a\sigma}^{\mathbb W_6}$. It is clear that $R \cup S$ is a subrack decomposition of  $\mathcal O_{a\sigma}^{\mathbb W_6}$. Consequently $\mathcal O_{a\sigma }^{\mathbb W_6}$ is of type {\rm D}. \hfill $\Box$

\begin {Lemma}\label {1.3''} If $\sigma$ is of type ($2^2, 3^1$),  then $\mathcal O_{a\sigma}^{\mathbb W_7}$ is of type ${\rm D}$ for all $a \in \mathbb Z_2^7.$
\end {Lemma}
\noindent {\it Proof.} Let $\pi=(5\ 6\ 7)$,  $\xi=(1\ 2)(3\ 4)$,   $\lambda=(1\ 3)(2\ 4)$,  $\tau = \pi \xi$ and $\mu=\pi \lambda$.
If $a=(a_1, a_2, a_3, a_4, 0, 0, 0)$,  let $b=(b_1, b_2, b_3, b_4, 1, 1, 0)$. If $a=(a_1, a_2, a_3, a_4, 1, 1, 1)$,  let $b=(b_1, b_2, b_3, b_4, 1, 0, 0)$. Then
the $5$-th,  $6$-th and $7$-th components of (\ref {e2.1.3}) are $(a_6+a_7, a_7+a_5, a_5+a_6)=(b_6+b_7, b_7+b_5, b_5+b_6)$, respectively. Consequently,  (\ref {e2.1.3}) does not hold.

Let $R := \mathbb Z_2^7 \rtimes \tau \cap \mathcal O_{a\sigma}^{\mathbb W_7}$ and $S := \mathbb Z_2^7 \rtimes \mu \cap \mathcal O_{a\sigma}^{\mathbb W_7}$. It is clear that $R \cup S$ is a subrack decomposition of  $\mathcal O_{a\sigma }^{\mathbb W_7}$. Consequently $\mathcal O_{a\sigma }^{\mathbb W_7}$ is of type $D.$ \hfill $\Box$

\begin {Example}\label {10.3}   If $\tau$ and $\mu$ are  of type $( 2^2)$, then $a\tau$ and $b\mu$ are  square commutative when they are conjugate to each other or $(-1) ^ {\sum _i^4 a_i} = (-1) ^ {\sum _i^4 b_i} $.

\end {Example}
\noindent {\it Proof.}    Let $\tau = (1\ 2) (3\ 4)$. Then  $\mu = (1\ 3) (2\ 4)$  or $\mu = \tau.$ It is clear that equation (\ref {e2.1.3}) becomes $a_1 + a_2 + a_3 +a_4
= b_1 +b_2 + b_3 + b_4.$  \hfill $\Box$

 \begin {Lemma}\label {2.6} {\rm (i)} Assume that $\tau$ and $\mu$ are conjugate  with ${\rm sq}(\tau, \mu) \not=\mu$ in $\mathbb S_n$ and $\tau(n) = \mu(n) =n$. If  $a \in \mathbb Z_2^n$ and there exists $i$ such that $a_i\not=a_n$ with $\tau(i)=\mu (i)= i$, then $\mathcal  O_{a\tau}$ is of type $D.$

{\rm (ii)} Assume that $n>4$ and $a\tau\in\mathbb W_n$ with type $(1^{n-2}, 2)$ of $\tau.$  If there exist $i, j$ such that $\tau(i)=i$ and $\tau(j)=j$ with $a_i\not=a_j$,  then $\mathcal O_{a\tau}$ is of type {\rm D}.

{\rm (iii)} Assume that $n>5$ and $a\tau\in\mathbb W_n$ with type $(1^{n-3}, 3)$ of $\tau.$ If there exist $i, j$ such that $\tau(i)=i$ and $\tau(j)=j$ with $a _i\not=a_j$,  then $\mathcal O_{a\tau}$ is of type {\rm D}.

\end {Lemma}
 \noindent {\it Proof.} {\rm  (i)}
 Let $R:= \{d\xi \in \mathcal O _{a \tau}^{\mathbb W _n}\mid \xi (n)=n; d_n=0\}$ and $S:= \{d\xi\in\mathcal O_{a\tau} ^{\mathbb W_n}\mid \xi (n) =n; d_n=1\}$. Obviously,  $R\cup S$ is a subrack decomposition.
 Let $\xi \in \mathbb S_n$ with $\xi (i) = i$ and $\xi (n) = n$ such that $\xi \rhd \tau  = \mu.$ By simple computation we have $(i, n) \xi \rhd (a \tau ) = b \mu$ with $b_n = a_i.$ Consequently, $a\tau$ and $b\mu$ are in the same set of  $R$ and $S$, which implies that $R\cup S$ is of type ${\rm D}$.

{\rm (ii)} It is clear that ${\rm sq}(\tau, \mu)\not=\mu$ with  $\tau :=(1, 2)$ and $\mu :=(2, 3)$. Applying Part (i) we complete the proof.

{\rm (iii)} It is clear that ${\rm sq}(\tau, \mu)\not=\mu.$ with $\tau=(1, 2, 3)$ and $\mu=(2, 4, 3)$. Applying Part (i) we complete the proof.
 \hfill $\Box$

\begin {Lemma}\label {1.2} Let $\sigma=(a, \tau)\in \mathbb W_n$.  If  $\mathcal{O}_{\tau}^{\mathbb S_n}$ is of type {\rm D},  then so is $\mathcal{O}_{(a, \tau)}^{\mathbb W_n}$.
\end {Lemma}
\noindent {\it Proof.} Let $X=S\cup T$ be a subrack decomposition of $\mathcal{O}_{\tau}^{\mathbb S_n}$ and $s\in S, t\in T$ such that
\begin {eqnarray} \label {e1.1}s\rhd (t \rhd (s \rhd t)) \not= t.\end {eqnarray}
Let $h \rhd \tau=s$ and $g\rhd \tau=t$ with $h, g\in {\mathbb S_n}$. It is clear ${\rm sq }((h\cdot a,  s), (g\cdot a,  t))\not=(g \cdot a,  t) $
since (\ref{e1.1}); $(h\cdot a, s)=h\rhd(a, \tau )$ and  $(g\cdot a, t)=g\rhd(a, \tau)$.

$(<{\mathbb S_n}\cdot a>,  X)$ is a subrack,  where $<{\mathbb S_n}\cdot a>$  is the subgroup generated by subset ${\mathbb S_n}\cdot a$ of $\mathbb Z_2^n$. In fact,  for any $h,  g\in {\mathbb S_n},  \xi,  \mu \in X, $  we have
\begin {eqnarray*}(h\cdot a, \xi)\rhd (g\cdot a, \mu)=((h +\xi g+\xi\mu\xi^{-1}h) \cdot a, \xi\rhd\mu )\in (<{\mathbb S_n}\cdot a>,  X).\end {eqnarray*}
Thus $(<{\mathbb S_n}\cdot a>,  X)$ is a subrack. Consequently  $(<{\mathbb S_n}\cdot a>, X) \cap \mathcal{O}_{(a,  \tau)}^{\mathbb W_n}=(<{\mathbb S_n}\cdot a>, S)
\cap \mathcal{O}_{(a,  \tau)}^{\mathbb W_n}\cup (<{\mathbb S_n}\cdot a>,  T) \cap\mathcal{O}_{(a,  \tau)}^{\mathbb W_n}$ is a subrack decomposition of $\mathcal{O}_{(a,  \tau)}^{\mathbb W_n}$ and is of type {\rm D}. \hfill $\Box$

\begin{Theorem}\label{1.4}
Let $G = \mathbb W_n$ with $n>4$. Let $\tau\in\mathbb S_n$ be of type $(1^{\lambda_1}, 2^{\lambda_2},  \dots,
n^{\lambda_n})$ and $a\in \mathbb Z_2^n$ with $\sigma=(a, \tau)\in G$ and $\tau\not=1$. If $\mathcal{O}_{\sigma}^G$ is not of
type {\rm D},  then the type of $\tau$ belongs to one in the following list.
\renewcommand{\theenumi}{\roman{enumi}}   \renewcommand{\labelenumi}{(\theenumi)}
\begin{enumerate}
\item $(2,  3)$; $(2^3);$
\item   $(2^4);$ $(1,  2^2), (1^2, 3), (1^2, 2^2) ;$
\item     $(1^{n-2},  2)$ and  $(1^{n-3},  3)$ $(n >5)$ with $a_i = a_j$ when $\tau (i) = i $ and $\tau (j)=j$.
\end{enumerate}
\end{Theorem}
\noindent {\it Proof.} It follows from Lemma \ref {1.3''}, Lemma \ref {2.6}, Lemma \ref {1.2}, Lemma \ref {1.3}, Lemma \ref {1.3'} and \cite [Theorem 4.1] {AFGV08}. \hfill $\Box$

\section{Nichols algebras  of irreducible {\rm YD} module over $\mathbb W_n$}\label {nichols algebra}
In this section we show that  except in three  cases Nichols algebras of irreducible {\rm YD} modules over classical Weyl groups
$\mathbb W_n$ are infinite dimensional.

 Let ${\rm supp} M := \{g \in G \mid M_g \not= 0\}$ for $G$-comodule $(M, \delta )$, where $M_g := \{x \in M \mid \delta (x) = g \otimes x \}$.

We shall use the following facts:

\begin{Theorem} \label {}(\cite [Cor. 8.4] {HS08}) Let $n\in \mathbb N $, $n\ge 3$, and assume that $G=\mathbb{S}_n$ is the
  symmetric group.
  Let $U$ be a {\rm YD} module  over $G$. If $\mathfrak B (U)$ is  finite dimensional, then $U$ is an irreducible
  {\rm YD} module over $G$.
\end {Theorem}

\begin{Theorem}\label{th:sm-intro} ( \cite [Th. 1.1] {AFGV08})
    Let $m\ge 5$.  Let $\sigma\in \mathbb S_m$ be of type $(1^{n_1},2^{n_2},\dots,m^{n_m})$,
    let $  \mathcal O_\sigma $ be the conjugacy class of $\sigma$ and let $\rho=(\rho,V) \in
    \widehat{{S_m } ^{\sigma}}
    $. If $\dim \mathfrak B ( \mathcal O_\sigma, \rho) < \infty$, then
    the type of $\sigma$ and $\rho$ are in the following list:
    \renewcommand{\theenumi}{\roman{enumi}}\renewcommand{\labelenumi}{(\theenumi)}
    \begin{enumerate}
        \item\label{it:caso2}
            $(1^{n_1}, 2)$, $\rho_1 = {\rm sgn}$ or $\epsilon$, $\rho_2 ={\rm sgn}$.
        \item
            $(2, 3)$ in $\mathbb S_5$, $\rho_2 ={\rm sgn}$, $\rho_3= \overrightarrow{\chi_{0}}$.
        \item\label{it:caso222}
            $(2^3)$ in $\mathbb S_6$, $\rho_2=\overrightarrow{\chi_{1}}\otimes\epsilon$ or
            $\overrightarrow{\chi_{1}}\otimes {\rm sgn}$.
    \end{enumerate}
\end{Theorem}

\begin {Theorem}\label {1.5}

If $G$ is a finite group and  $\mathcal{O}_{\sigma}^G$ is of type {\rm D},  then {\rm dim} $\mathfrak B(\mathcal{O}_{\sigma}^G,  \rho) = \infty $ for any $\rho \in \widehat { G^\sigma}$.
\end {Theorem}
\noindent {\it Proof.} It follows from \cite [Theorem 3.6]{AFGV08}.  \hfill $\Box$

\begin{Theorem}\label{1.6} Assume $n>4$. Let $\tau\in\mathbb S_n$ be of type $(1^{\lambda_1},  2^{\lambda_2}, \dots, n^{\lambda_n})$ and $a\in\mathbb Z_2^n$ with $\sigma=(a, \tau)\in \mathbb W_n$ and $\tau\not=1$. If {\rm dim} $\mathfrak B(\mathcal{O}_{\sigma}^{\mathbb W_n}, \rho)<\infty $,  then the type of $\tau$ belongs to one in the following list:
\renewcommand{\theenumi}{\roman{enumi}}   \renewcommand{\labelenumi}{(\theenumi)}
\begin{enumerate}
\item $(2,  3)$; $(2^3);$
\item   $(2^4);$ $(1,  2^2), (1^2, 3), (1^2, 2^2) ;$
\item     $(1^{n-2},  2)$ and  $(1^{n-3},  3)$ $(n >5)$ with $a_i = a_j$ when $\tau (i) = i $ and $\tau (j)=j$.
\end{enumerate}
\end{Theorem}
\noindent {\it Proof.} It follows from Theorem \ref {1.4} and Theorem \ref {1.5}. \hfill $\Box$

\vskip.1in
 Let $a\mu =c\tau \# d\xi \in\mathbb W_n$ with $c\tau \in \mathbb W_m$,  $d\xi \in \mathbb
B_{n-m}$ and   $c\tau \bot d\xi$. Let
  $\rho=\rho _1\otimes \rho_2 \in  \widehat {
\mathbb W_n ^{a\mu}}=\widehat {\mathbb W_m ^{c\tau}} \times\widehat {\mathbb W_{n-m} ^{d\xi}}$.
If  $c=(1, 1, \cdots, 1)$, then    $\rho_1=(\chi _1\otimes\mu _1)\uparrow _{
G^{a\mu}_{\chi_1}} ^{G^{a\mu}}$ with $\chi_1 \in \widehat {(Z_2^ {m})^\tau} $,  $\mu _1\in \widehat{(\mathbb S_m^\tau)_{\chi_1}}$ (see \cite [provious Proposition 2.9 in Section 2.5] {ZZ12} ), where $G= \mathbb W_n $.
Case $a=0$ and $a=(1, 1, \cdots, 1)$ were studied in
paper \cite [Theorem 1.1 and Table 1]{ZZ12}. Other cases are listed  as follows:

\begin {Corollary}\label {2.9} Under notation above assume  $\dim \mathfrak B( {\mathcal O} _{a\mu} ^{\mathbb W_n},  \rho)< \infty$ with  $d_1=d_2=\cdots=d_{n-m}$ and   $\xi ={\rm id}^{\otimes (n-m)}$.

{\rm (i)} Then $\rho _1 (c\tau)=-id $ when $\rho _2(d\xi)={\rm id} $ and $\rho _1(c\tau)={\rm id}$ when $\rho _2(d\xi)=-{\rm id} $.

{\rm (ii)} Case $\tau=(1\ 2)$,  $c=0$ and $d=(1, 1, \cdots, 1)$. Then $\rho_1(c\tau)=\pm {\rm id}$,  $\rho_2(d\xi)=\mp{\rm id} $,  $\chi_1(c)=1$ and $\mu _1 (\tau)=\pm 1.$

{\rm (iii)} Case $\tau=(1\ 2)$,  $c=(1, 1)$ and $d =0$. Then  $\rho_1 (c\tau) =-{\rm id}$,  $\rho _2 (d\xi)={\rm id}$.

{\rm (iv) } Case $\tau=(1\ 2\ 3)$,  $c=0$ and  $d =(1, 1, \cdots, 1)$. Then $\rho_1(c\tau)={\rm id}$,  $\rho_2(d\xi)=-id$,  $\chi_1(c)=1$ and $\mu _1 (\tau)= 1.$

{\rm (v)} Case $\tau=(1\ 2\ 3)$,  $c=(1, 1, 1)$ and $d =0$. Then $\rho_1(c\tau)=-id$,  $\rho_2(d\xi)={\rm id}$,  $\chi_1(c)=-1$ and $\mu _1 (\tau)=1.$

{\rm (vi)} Case $\tau=(1\ 2)( 3\ 4)$,  $c=0$ and $d=(1, 1)$. Then $\rho_1(c\tau)={\rm id}$,  $\rho_2(d\xi)=-id$.

{\rm (vii)} Case $\tau=(1\ 2)( 3\ 4)$,  $c=(1, 0, 1, 0)$ and $d=(1, 1)$. Then $\rho_1(c\tau)=\pm {\rm id}$,  $\rho_2(d\xi)=\mp {\rm id}$.

{\rm (viii)} Case $\tau=(1\ 2)( 3\ 4)$,  $c=(1, 0, 1, 0)$ and $d=0$. Then $\rho_1(c\tau)=-id$,  $\rho _2(d\xi)={\rm id}$.

{\rm (ix)} Case $\tau=(1\ 2)( 3\ 4)$,  $c=(1, 0, 0, 0)$ and $d=(1, 1)$. Then $\rho_1(c\tau)=\pm {\rm id} $,  $\rho_2(d\xi)=\mp {\rm id}$.

{\rm (x)} Case $\tau =(1\ 2)( 3\ 4)$,  $c=(1, 0, 0, 0)$ and $d=(0, 0)$. Then $\rho_1(c\tau)=-id$,  $\rho_2(d\xi)={\rm id}$.
\end {Corollary}
\noindent {\it Proof.} {\rm (iv)} If $\rho_2(d\xi)={\rm id}$,  then $\dim \mathfrak{B}({\mathcal O}_{c\tau}^{\mathbb W_{3}}, \rho_1)<\infty$ by Lemma \ref {2.7},  which constracts to \cite [Theorem 1.1] {ZZ12}.

{\rm (v)} If $\chi _1(c) =1$,  then there exists a contradiction by \cite [Proposition 2.4,  Theorem 1.1] {ZZ12}.

The others  follow from Lemma \ref {2.7} and \cite [Theorem 1.1]{ZZ12}. \hfill $\Box$

\appendix

\section*{Appendix}
In this section we give the general FK' Conjecture and another proof of Theorem \ref {1.5}.

\section {General FK' Conjecture.}

$\mathfrak{B}({\mathcal O}_{\sigma}, \rho)$ is finite dimensional when $n\le 5$ according to \cite {MS,  FK99,  AZ07}. However,  it is open whether $\mathfrak{B}({\mathcal O}_{\sigma}, \rho)$ is finite dimensional when $n>5.$
We know the relationship between  FK' Conjecture  and Nichols algebra
$\mathfrak{B} ({\mathcal O}_{{(1, 2)}} ,  \epsilon \otimes {\rm sgn})$ of transposition over symmetry group(see \cite
[Theorem 5.7]{MS}). That is,  if  $\dim \mathfrak{B} ({\mathcal O}_{{(1, 2)}} , \epsilon\otimes {\rm sgn}) = \infty$,  then so is  $\dim \mathcal E_n$.

\begin {Definition} \label {4.4'}  (See \cite [Def. 2.1]{FK99}) algebra ${\mathcal E}_n$ is generated by $\{x_{ij} \mid 1\le i< j \le n\}$ with definition relations:

{\rm (i)} $x_{ij}^2=0$ for $i<j.$

{\rm (ii)}  $ x_{ij}x_{jk}=x_{jk}x_{i k}+x_{i k}x_{ij}$ and $x_{jk}x_{ij}=x_{i k}x_{jk}+x_{ij}x_{i k}, $ for $i<j<k.$

{\rm (iii)} $x_{ij}x_{kl}=x_{kl}x_{ij}$  for any distinct $i, j, k, l$ and $l, $ $i<j, k<l.$
\noindent
Equivalently,  algebra ${\mathcal E}_n$ is generated by $\{x_{ij}\mid i\not=j, 1\le i, j \le n\}$ with definition relations:

{\rm (i)} $x_{ij}^2=0$,  $x_{i j}=-x_{ji}$;

{\rm (ii)} $ x _{ij}x_{jk}+x_{jk}x_{ki}+x_{ki}x_{ij}=0;$

{\rm (iii)} $x_{ij}x_{kl}=x_{kl}x_{ij}$;

for any distinct $i, j, k$ and $l.$
\end {Definition}

S. Fomin and A.N. Kirillov \cite [Conjecture 2.2]{FK99} point out a conjecture `` ${\mathcal E}_n$ is finite dimensional ".   Consequently,    we have

\begin {Problem} \label {4.4} Let $\alpha_{ijk}, \beta_{ijk}, \gamma_{ij}, \lambda_{ijkl}\in\{1, -1\}$ for any distinct $i, j, k$ and $l$ with $1\le i,  j k,  l \le n$. Assume that algebra $A(\alpha, \beta, \gamma, \lambda )$ is generated by $\{x_{ij}\mid i\not=j, 1\le i, j\le n\}$ with definition relations:

  {\rm  (i)} $x_{ij}^2=0$,  $x_{i j}=\gamma_{ij}x_{ji}$;

 {\rm  (ii)} $x_{ij}x_{jk}+\alpha_{ijk}x_{jk}x_{ki}+\beta_{ijk}x_{ki}x_{ij}=0;$

  {\rm (iii)} $x_{ij}x_{kl}=\lambda_{ijkl}x_{kl}x_{ij}$;

for any distinct $i, j, k$ and $l.$

Then $A(\alpha, \beta, \gamma, \lambda )$ is finite dimensional when $n>5$.

\end {Problem}
Obviously,  $\mathcal E_n=A(1, 1, -1, 1), $ i.e. $\mathcal E_n= A(\alpha, \beta, \gamma, \lambda )$ with $\alpha _{i,  j,  k} = 1, \beta _{i,  j,  k} =1, \gamma _{ij} =-1, \lambda _{i,  j,  k,  l} =1$ for any distinct $i, j, k$ and $l.$

\section {Proof of Theorem \ref {1.5}.}

 \noindent {\it Proof.} {\rm (i)}  There exists subrack $X \subseteq \mathcal{O}_{\sigma}^G$.
$R \cup S= X$ is a subrack  decomposition of $X$there exists $a\in R$, $b\in S$ such that ${\rm sq} (a, b) := a\rhd (b \rhd (a \rhd b)) \not= b$. Let $H$ be the subgroup generated by $X$ in  $G$. It is clear ${\mathcal O}_{\sigma _1}^H \not= {\mathcal O}_{\sigma _2} ^H$ for any $\sigma _1 \in R, \sigma _2 \in S. $ In fact, if ${\mathcal O}_{\sigma _1}^H = {\mathcal O}_{\sigma _2} ^H$, then there exists $t\in H$, such that $\sigma _2 = t \sigma _1t^{-1}$. $t = \prod _{i=1} ^n t_i$ whit $t_i \in X$ or $t_i^{-1}\in X.$ Therefore $\sigma_2 \in R$, which is a contradiction.

{\rm (ii)} By {\rm (i)}, we have ${\mathcal O}_{a}^H \not= {\mathcal O}_{b} ^H$  and they are not square commutative.  By \cite [Th. 8.6] {HS08}, $\dim \mathfrak B (M({\mathcal O}_{a}^H, \lambda  )+ M({\mathcal O}_{b} ^H), \xi )) = \infty,  $ for any
$\lambda \in \widehat {H^a}$, $\xi  \in \widehat {H^b}$. Consequently,  {\rm dim} $\mathfrak B(\mathcal{O}_{\sigma}^G,  \rho) = \infty $ for any $\rho \in \widehat { G^\sigma}$ by  \cite [Lemma 3.2{ \rm ii}]{AFGV08}. \hfill $\Box$

\end {document}